\documentclass[a4paper]{article}

\usepackage{graphicx}
\usepackage{amssymb}
\usepackage{latexsym}
\usepackage{verbatim}
\usepackage{amsthm}

\theoremstyle{definition}
\newtheorem{example}{Example}[section]
\newtheorem{definition}{Definition}[section]
\theoremstyle{plain}
\newtheorem{theorem}{Theorem}[section]

\newtheorem{lemma}{Lemma}[section]

\def\sqr#1#2{{\vcenter{\vbox{\hrule height.#2pt
\hbox{\vrule width.#2pt height#1pt \kern#1pt \vrule width.#2pt}
\hrule height.#2pt}}}}
\def\square{\mathchoice\sqr55\sqr55\sqr{2.1}3\sqr{1.5}3}

\title{\bf Duality on hypermaps with symmetric or alternating monodromy group}

\author{Daniel Pinto\\
\small CMUC, Department of Mathematics, University of Coimbra \\[-0.8ex]
\small 3001-454 Coimbra, Portugal \\
\small \texttt{dpinto@mat.uc.pt} }

\begin{document}

\maketitle

\begin{abstract}

Duality is the operation that interchanges hypervertices and
hyperfaces on oriented hypermaps. The duality index measures how far a
hypermap is from being self-dual. We say that an oriented regular hypermap has \emph{duality-type}
$\{l,n\}$ if $l$ is the valency of its vertices and $n$ is the
valency of its faces. Here, we study some properties
of this duality index in oriented regular hypermaps and we prove that for each pair $n$, $l \in \mathbb{N}$, with $n,l \geq 2$, it is
possible to find an oriented regular hypermap with extreme duality
index and of duality-type $\{l,n \}$,  even if we are restricted to hypermaps with alternating or symmetric monodromy group.

\medskip 



\end{abstract}

\section{Introduction}

One can look at a hypergraph as a generalization of the well known
notion of a graph, since an important restriction is removed: in
a hypergraph, an (hyper)edge might connect more than two
(hyper)vertices. In a similar way, a hypermap is an obvious
generalization of the concept of a map. Although the word
\emph{map} is often used in mathematics with a different meaning,
here it is defined as a cellular embedding of a connected graph on
a compact connected surface. Without surprise, it is somehow
natural to say that a hypermap is, in its topological form, a
cellular embedding of a connected hypergraph. Or, alternatively,
an embedding of a connected trivalent graph on a compact connected
surface, labelling each face 0, 1 or 2, so that each edge of the
trivalent graph is incident to two faces carrying different
labels. In this representation of the hypermap (James
representation \cite{James}), the hypervertices are usually
represented by label 0, the hyperedges by label 1, and the
hyperfaces by label 2. We should note that, by doing this, each
hyperedge of the hypermap might be adjacent to more than two
hypervertices and, consequently, have valency greater than 2
(something that does not occur in a map). If the valencies of all
hyperfaces are equal, and the same happens to the hypervertices
and to the hyperedges, we say that the hypermap is \emph{uniform}.
Here, all the surfaces will be orientable and without boundary. It
follows that every hypermap can give rise to two \emph{oriented}
hypermaps, one for each fixed orientation.  If there is an
automorphism that sends one of those oriented hypermaps into the
other, we say that the hypermap is \emph{reflexible}. Otherwise,
we call it \emph{chiral}. The topological idea of an oriented
hypermap (a cellular embedding of a hypergraph on an oriented
 surface), briefly sketched here, has a very useful algebraic
translation. In fact, one can look at the results of this paper as
pure algebraic results, although they certainly have an important
topological meaning.

\medskip

A hypermap can be regarded as a transitive permutation
representation $\Delta\to{\rm Sym}\,\Omega$ of the group
$$\Delta=\langle r_0, r_1, r_2\mid r_0^2=r_1^2=r_2^2=1\rangle\cong
C_2*C_2*C_2,$$ on a set $\Omega$ representing its hyperflags (the
cells of the barycentric subdivision of the hypermap). Similarly,
an oriented hypermap (without boundary) can be regarded as a
transitive permutation representation of the subgroup
$$\Delta^+=\langle \rho_0, \rho_1, \rho_2\mid \rho_0\rho_1\rho_2
=1\rangle=\langle \rho_0,\rho_2\mid -\rangle$$ of index $2$ in
$\Delta$ (a free group of rank $2$) consisting of the elements of
even word-length in the generators $r_i$, where $\rho_0=r_1r_2$,
$\rho_1=r_2r_0$ and $\rho_2=r_0r_1$. In the case of hypermaps, the
hypervertices, hyperedges and hyperfaces ($i$-dimensional
constituents for $i=0, 1, 2$) are the orbits of the dihedral
subgroups $\langle r_1, r_2\rangle$, $\langle r_2, r_0\rangle$ and
$\langle r_0, r_1\rangle$, and in the case of oriented hypermaps
they are the orbits of the cyclic subgroups
$\langle\rho_0\rangle$, $\langle\rho_1\rangle$ and
$\langle\rho_2\rangle$, with incidence given by nonempty
intersection in each case. The local orientation around each
hypervertex, hyperedge or hyperface is determined by the cyclic
order of the corresponding cycle of $\rho_0, \rho_1$ or $\rho_2$.

\medskip

An oriented hypermap is regular if it has, quoting J. Siran
\cite{Siran}, the \emph{highest degree of symmetry}. It is well
known that, algebraically, an oriented regular hypermap $\cal H$
can be represented by a triple $(G,x,y)$, where the group $G$, a
quotient of $\Delta^+$, is generated by the elements $x$ and $y$.
If we want to look at those generators from a topological point of
view, $x$ can be interpreted as the permutation that cyclic
permutes the hyperdarts (oriented hyperedges) based on the same
hypervertex, and $y$ the permutation that cyclic permutes the
hyperdarts based on the same hyperface, according to the chosen
orientation. It is not difficult to verify that every regular
hypermap is uniform (in the James representation: every face
labelled $i$, with $i \in \{0,1,2\}$, has the same valency). We
say that a regular hypermap has \emph{type} $(l,m,n)$ (for $l$,
$m$, $n \in \mathbb{N}$) if $l$, $m$ and $n$ are the valencies of
its hypervertices, hyperedges and hyperfaces, respectively. The
group $G=\langle x,y \rangle$ is called the \emph{monodromy group}
of the hypermap (and, since the hypermap is regular, it coincides
with its automorphism group). Duality, one of the many possible
operations on hypermaps (see \cite{GarethDaniel}), is the
operation that interchanges hyperfaces and hypervertices, which in
terms of the generators of the monodromy group is the same as
saying that the dual of ${\cal H}=(G,a,b)$ is ${\cal
H}^d=(G,b,a)$. If there is an automorphism of $G$ that
interchanges $a$ and $b$, we say that ${\cal H}$ is
\emph{self-dual} (invariant under duality) and ${\cal H} \cong
{\cal H}^d$.

\section{Duality index}

A possible way to define the \emph{duality group} of a hypermap
${\cal H}=(G,x,y)$, as an extension of the notion of chirality
group (see \cite{ChiralityGroup}, \cite{GarethDaniel}), is to say
that it is the minimal subgroup $D({\cal H}) \trianglelefteq G$
such that the quotient hypermap: $${\cal H}/D({\cal H})=
(G/D({\cal H}), xD({\cal H}), yD({\cal H}))$$ is a self-dual
hypermap. The order of this duality group is called the
\emph{duality index}\index{duality!index} $d$ of $\cal H$, and it
is a way to measure how far a hypermap is from being self-dual. If
the duality index is 1, then the hypermap is self-dual; and the
bigger that index, the more distant the hypermap is from being
self-dual. If $D({\cal H}) = G$, the duality group is equal to the
monodromy group of the hypermap and we say that the hypermap has
\emph{extreme duality index}.

\bigskip

\section{The symmetric and the alternating groups}

Before dealing with more difficult problems, we start this section
by proving a simple lemma:

\begin{lemma}
For every $n \in \mathbb{N}$:
\begin{itemize}
\item [a)] there is a non-self-dual hypermap with monodromy group
$S_n$.

\item [b)] there is a non-self-dual hypermap with monodromy group
$A_n$.

\end{itemize}
\end{lemma}
\noindent Proof: a) If $n>2$, we take $x=(1,2,\quad...\quad, n)$
and $y=(1,2)$. These permutations generate the group $S_n$ but,
because they have different orders, there is no automorphism that
interchanges those two generators. It follows that the hypermap
${\cal H}=(S_n,x,y)$ is not self-dual. If $n=2$, $S_2 \cong C_2=
\langle g \rangle$ and $(C_2,1,g)$ is not self-dual (if we take 1
as the neutral element of $C_2$).

\medskip

b) For $n > 3$, we just need to apply that same idea but using the
canonic generators $x=(1,2,...,n)$ and $y=(1,2,3)$, if $n$ is odd
and greater than 3, or $x=(2,3,..,n)$ and $y=(1,2,3)$ if $n$ is
even. Because, in both cases, the two generators do not have the
same order, the hypermap is not self-dual. If $n=3$, $A_n=A_3\cong
C_3=\langle g \rangle$ and $(C_3,1,g)$ is not self-dual. $\hfill
\square$

\medskip

\begin{example} Let $G=\langle x=(1,2,3,4,5),y=(1,2,3)\rangle =A_5$. Because $A_5$ is simple and the
hypermap $(A_5,x,y)$ is not self-dual, the duality group (which is
normal in $A_5$) must be $A_5$ itself. Hence, the duality index of
the hypermap is $|A_5|=5!/2$
\end{example}

More generally, if we take $$(A_n, (1,2,...,n),(1,2,3)), \mbox{ if
 }n \mbox{ is odd, or }
 A_n, (2,3,...,n), (1,2,3)), \mbox{ if } n \mbox{ is even, }$$ we
get hypermaps with extreme duality index and with monodromy group
$A_n$. Hence:

\begin{theorem}
If $n \geq 3$ there is a hypermap with extreme duality index and
with monodromy group $A_n$. $\hfill \square$
\end{theorem}

\noindent \emph{Remarks}: a) Since each one of those hypermaps
$\cal H$ has extreme duality index, $D({\cal H})= Mon ({\cal
H})=A_n$. It follows that, for any $n \geq 3$, $A_n$ is the
duality group of some hypermap. b) If we take any simple group
generated by two elements of different orders, we can get a
similar result.

\medskip

If $G=S_n$, the only possible duality groups are $1$, $A_n$ and
$S_n$. Therefore, a hypermap with monodromy group $S_n$ is
self-dual, has extreme duality index or has $n!/2$ as duality
index. For $n \neq 6$, all automorphisms of $S_n$ are inner and
act by conjugation. If there is an automorphism that transposes
the two generators,  the hypermap is self-dual and then ${\cal
D}({\cal H})=1$. Otherwise, ${\cal D}({\cal H})=A_n$ or $S_n$ and
we need to check if the hypermap with monodromy group $S_n/A_n$,
of order 2, is self-dual or not.

\bigskip

\begin{theorem}\label{sn}
Every hypermap ${\cal H}=(S_n,x,y)$:
\begin{itemize}

\item[i)] has extreme duality index if $x$ or $y$ is an even
permutation;

\item[ii)] is self-dual or has duality index $n!/2$ if $x$ and $y$
are both odd permutations and $ n \neq 4$.

\item[iii)] is self dual or has duality index 4 if $x$ and $y$ are
both odd permutations and $n=4$.

\end{itemize}
\end{theorem}
\noindent Proof: i) The only non-identity quotients $S_n/N$
of $S_n$ are $S_n/A_n \cong C_2$ and $S_4/V_4 \cong S_3$, when
$n=4$. In each case, because one of $xN$ and $yN$ is in the unique
subgroup of index 2 of $S_n/N$ and the other is not, there can be
no automorphism of $S_n/N$ transposing $xN$ and $yN$. So, the only
self-dual quotient is the trivial one and the hypermap has extreme
duality index.

\vskip 0,2cm

 ii) $S_n$ is not cyclic and, by definition $\langle
x,y\rangle =S_n$. Hence, $x \neq y$. Suppose $\cal H$ is not
self-dual. Because ${\cal H}/A_n=(S_n/A_n, xA_n, yA_n)$ and
$|S_n/A_n|=2$, the two generators $xA_n$ and $yA_n$ must be the
same. It follows that the hypermap ${\cal H}/A_n$ is self-dual and
$|D({\cal H})|=n!/2$.

\vskip 0,2cm

iii) Let $x$ and $y$ be two odd permutations generating $S_4$ and
$N$ the Klein group $V_4$ (a normal subgroup in $S_4$). That pair
of generators may be formed by a transposition and a 4-cycle or by
two 4-cycles. Say $x$ is a transposition and $y$ is a 4-cycle.
Then, they map to distinct permutations in $S_4/N \cong S_3$.
Because the third element in $S_3$ conjugates each to the other,
the quotient hypermap is self-dual (whereas the hypermap itself is
not), with duality group $V_4$. If $x$ and $y$ are both 4-cycles
the hypermap is self-dual. $\hfill \square$

\medskip
\textbf{Can we have an oriented regular hypermap with extreme
duality index for each type $(l,m,n)$? }

\medskip
The answer is no. There are no hypermaps of extreme duality index
and of type $(3,2,3)$, since this has to be the tetrahedron, which
is self-dual. On the other hand, if we restrict ourselves to
hyperbolic triples, where $l^{-1}+m^{-1}+n^{-1}<1$, we can
enumerate orientably regular hypermaps of a given type $(l,m,n)$
with automorphism groups isomorphic to $PSL(2,q)$ or $PGL(2,q)$
(this enumeration can be found in a joint work of Marston Conder,
Primoz Potocnik and Josef Siran \cite{projective}, based on a
paper by Sah \cite{Sah}). Because these groups are simple or
almost simple, we can use them to try to find hypermaps with
extreme duality index or self-dual hypermaps. If $l \neq n$ the
hypermap cannot be self-dual. If $l=n$, we have to check if there
is an automorphism of $PSL(2,q)$ or $PGL(2,q)$ that interchange
the two generators.

We should also notice that for some triples $(l,m,n)$ is very easy
to find hypermaps with extreme duality index and of that type.
\medskip

If $g.c.d.(l,n)=1$ then $$G_{l,m,n}=\langle
a,b|a^{l}=b^{n}=(ab)^{m}=1\rangle$$ is obviously the monodromy
group of a hypermap of type $(l,m,n)$. Let $N$ be the duality group of
the hypermap $(G_{l,m,n},a,b)$. Then, $N$ is the smallest normal
subgroup of $G_{l,m,n}$ such that $G_{l,m,n}/N$ is reflexible,
which means that we are working with the smallest normal subgroup
$N$ of $G_{l,m,n}$ such that the assignment $a\mapsto b$, $b
\mapsto a$ induces an automorphism of $G_{l,m,n}/N$. We obtain
this quotient by adding extra relations, substituting $a$ for $b$
and $b$ for $a$ in the original ones. Then:
$$ G_{l,m,n}/N=\langle
a,b|a^{l}=a^{n}=b^{n}=b^{l}=(ab)^{m}=(ba)^m=1\rangle.$$ Because
$l$ and $n$ are co-prime, the group $G_{l,m,n}/N$ collapses to the
identity. Therefore, ${\cal H}=(G_{l,m,n},a,b)$ is a hypermap with
extreme duality index and of type $(l,m,n)$. This proves that, if
$g.c.d.(l,n)=1$, for every triple $(l,m,n)$ there is a hypermap
with extreme duality index and of that type $(l,m,n)$.

\subsection{Alternating monodromy group}
In the 60's, Graham Higman conjectured that any Fuchsian group
has among its homomorphic images all but finitely many of the
alternating groups. He also proved that $A_n$ is a factor group of
$(2,3,7)=\langle a,b|a^2,b^3,(ab)^7 \rangle$ for all large $n$.
Because 2, 3 and 7 are prime numbers and $2 \ne 7$, we can
conclude, from that result, that there is an infinite number of
 hypermaps with extreme duality index and of type $(2,3,7)$. The result obtained by
Higman was later extended by others (Marston Conder
\cite{Conder2}, for instance) that proved that the same can be
said for other families of triangle groups. The complete proof of
the conjecture, nevertheless, was only published in 2007 by Brent
Everitt \cite{Everitt}. In his paper, it is shown that we only
need to consider the triangle groups $(p,q,r)$, $3 \leq p <q<r$ to
prove the main result (which is done making use of coset diagrams
for those triangle groups). Hence, it is possible to say that if
$l$, $m$, $n$ are prime and $l \neq n$ we can always find infinite
hypermaps with extreme duality index and of type $(l,m,n)$, with
alternating group as monodromy group. If $p$, $q$, $r$ are not all
prime, the alternating groups, being factor groups of the triangle
group $(p,q,r)$, might correspond to hypermaps of type
$(p',q',r')$ with $p'|p$, $q'|q$ and $r'|r$, and not always of
type $(p,q,r)$. It follows that we cannot directly apply the
result of Brent Everitt to prove that there are hypermaps of any
type with alternating monodromy group. The proof of Brent Everitt
would had to be changed, in order to sustain that conclusion.
However, here we will use a different approach.

\subsection{Duality-type of hypermaps with symmetric or alternating monodromy group}

\begin{definition}
We say that an oriented regular hypermap has \emph{duality-type}
$\{l,n\}$ if $l$ is the valency of its vertices and $n$ is the
valency of its faces (which is the same as saying that $l$ and $n$
are the orders of the generators interchanged by the duality
operation).
\end{definition}

If we restrict ourselves to the family of hypermaps whose
monodromy group is the alternating or the symmetric group, can we
have a hypermap with extreme duality index and of any duality-type
$\{l,n\}$? In fact, we can not only prove that the answer is
affirmative but also explicitly show how to construct those
hypermaps. The reason why, here, we look for \emph{duality-type}
instead of \emph{type} is because, in the first case, we just need
to control the order of the two generators $x$ and $y$, while in
the second case we also need to pay attention to the order of
$xy$, which seems a harder problem to solve. To complete that
easier first task, we still need to add a few concepts and
theorems about primitive groups.

Let $(G,\Omega)$ be a permutation group. An equivalence relation
$\sim$ is called \emph{G-invariant} if whenever $\alpha, \beta \in
\Omega$ satisfy $\alpha \sim \beta$ then $g(\alpha) \sim g(\beta)$
for all $\alpha, \beta \in \Omega$. Two obvious $G$-invariant
equivalence relations are: (i) $\alpha \sim \beta$ if and only if
$\alpha=\beta$ and (ii) $\alpha \sim \beta$  for all $\alpha,
\beta \in \Omega$. We call $(G,\Omega)$ \emph{imprimitive} if it
admits some equivalence relation other than (i) or (ii).
Otherwise, we call $(G,\Omega)$ \emph{primitive}.

There are several examples of primitive groups (Colva M.
Roney-Dougal \cite{Colva} has classified all the  primitive
permutation groups of degree less than 2500). For instance, any
alternating group $A_n$ is primitive, and so is $PGL(2,q)$, in its
standard action, for any prime power $q$. However, it is known -
as underlined by Peter Cameron in the Encyclopedia of Design
Theory \cite{Encyclopaedia} - that primitive groups are ``rare (for
almost all $n$, the only primitive groups of degree $n$ are the
symmetric and alternating groups, see \cite{Cameron}); and small
(of order at most $n^{c\cdot log n}$ with known exceptions)''. The
fact that most of the primitive groups are alternating groups
(which are simple) and symmetric groups (which have only one non
trivial normal subgroup) makes \emph{primitivity} a powerful
concept to build hypermaps with extreme duality index. Although
the probability of a primitive group being an alternating or a
symmetric group is very high, we need to be sure that we are not
dealing with a different kind of group. The next definitions and
theorems help us to achieve that goal.
\bigskip
\bigskip

\begin{definition}
Let $G$ be a permutation group on $\Omega$ and $k$ a natural
number with $1 \leq k \leq n=|\Omega|$. $G$ is called
\emph{$k$-transitive} if for every two ordered $k$-tuples
$\alpha_1,...\alpha_k$ and $\beta_1,...,\beta_k$ of points of
$\Omega$ (with $\alpha_i \neq \alpha_j$, $\beta_i \neq \beta_j$,
for $i \neq j$) there exists $g \in G$ which takes $\alpha_t$ into
$\beta_t$ (for $t=1,...,k$).
\end{definition}

\begin{definition}
The \emph{degree} of a permutation group is the number of points
moved by at least one of its permutations.
\end{definition}

\begin{definition}
The \emph{minimal degree} of a permutation group is the minimum
number of points moved by any non-identity element of the group.
\end{definition}

\begin{theorem} \emph{\cite{Miller2}} \label{theorem a}   The minimal degree of a primitive group which is
neither alternating nor symmetric must exceed 4 whenever its
degree exceeds 8.
\end{theorem}

\begin{theorem} \emph{\cite{Wielandt}} \label{theorem b}A primitive group of degree $n$, which
contains a cycle of degree $m$ with $1<m<n$ is
$(n-m+1)$-transitive.
\end{theorem}

\begin{theorem} \emph{\cite{Rotman}} \label{theorem c}
A 2-transitive permutation group is transitive.
\end{theorem}

And, as a consequence of the classification of simple groups
\cite{Cameron2}:

\begin{theorem} \label{theorem d} If a permutation group $G$ is at least $6$-transitive then $G$
is the alternating group or the symmetric group.
\end{theorem}

\bigskip

\bigskip
\bigskip

With these tools we can now prove the following theorem:
\begin{theorem}
For each pair $n$, $l \in \mathbb{N}$, with $n,l \geq 2$, it is
possible to find an oriented regular hypermap with extreme duality
index and of duality-type $\{l,n \}$, with alternating or
symmetric group as monodromy group.
\end{theorem}

\noindent Proof: The general idea behind this proof is the
following: first, we try to prove that the two chosen permutations
generate a primitive group, and then we show that they must
generate the symmetric or the alternating group. If they generate
the symmetric group (which is not a simple group), we also need to
check that they give rise to a hypermap with extreme duality
index.

\begin{itemize}
\item[a)] If $l$, $n$ are both odd and $l \neq n$ (we may assume,
without loss of generality, that $l>n$) let $x=(1,2,...,l)$ and
$y=(1,2,...,n)$ be two permutations of $S_l$.

The only equivalence relations preserved by $x$ are the
congruences $mod \mbox{ } k_1$ (for some $k_1|l$). If that
equivalence relations are also preserved by $y$ then they have to
be the congruences $mod \mbox{ } k_2$ (for some $k_2|n$). Hence,
the only equivalence relations preserved by $x$ and $y$ are the
congruences $mod \mbox{ } k$ (for some $k|l$ and $k|n$). Then,
$$ l \equiv n (mod \mbox{ } k) \Rightarrow l \sim n.$$

If $z$ is the commutator $y^{-1}x^{-1}yx=(n,n-1,l)$ we have:
$$z(l) \sim z(n) \Leftrightarrow n \sim n-1$$
Therefore, because $(n-1,n)=1$, the only equivalence relation
preserved by $G=<x,y>$ is a trivial one.

Hence, $G$ is primitive. But $z \in G$ and the minimal degree of a
primitive group which is neither alternating or symmetric must
exceed 4 whenever its degree exceeds 8 (Theorem \ref{theorem a}).
Therefore, if $l>8$, $G$ must be the alternating group and ${\cal
H}=(G=A_l,x,y)$ is a hypermap of extreme duality index and of
duality type $\{l,n \}$. (For $l<8$ we can easily find two
permutations, one of order $l$ and another of order $n$, that
generate $A_l$.)

\item[b)] If $l$ is even and $n$ is odd (or if $l$ is odd and $n$
is even) and $l \neq n$ a similar proof can be written using the
cycles $x=(1,2,...,l)$ and $y=(1,2,...,n)$. In this case, however,
$G=<x,y>=S_l$. Because not both of the generators are odd
permutations, ${\cal H}=(G,x,y)$ has extreme duality index (see
Theorem \ref{sn}).

\item[c)] If $l$ and $n$ are both even and $l \neq n$ (assuming,
without loss of generality, that $l > n$), we take
$$x=(1,2,...,l)$$
$$y=(1,2)(l,l+1,l+2,...,l+n-1)$$
permutations of $S_{l+n-1}$ and $G=<x,y>$. First, we need to show
that $G$ is primitive. Let $\sim$ be a $G$-invariant equivalence
relation and suppose $B_1$ is an equivalence class (or a group
block) for $\sim$ such that $(l+1) \in B_1$.

\begin{itemize}
\item[i)] If $l \in B_1$ then, because $x(B_1)=B_1$ (since $l+1$
is fixed by $x$), all elements of the set $\{1,2,..,l\}$ must
belong to $B_1$.

Then, $B_1y=B_1$ (since $y$ fixes $3 \in B_1$). Therefore $(l+2),
(l+3),...,(l+n-1)$ also belong to $B_1$. It follows that all
elements of $\Omega=\{1,2,...,l+n-1\}$ belong to the same block
and that the equivalence relation must be trivial.

\item[ii)] If $l$ does not belong to $B_1$, we assume that $l \in
B_2 \neq B_1$.

Suppose $(l+2) \in B_2$. Then: $$x(l+2)=l+2$$
$$x(l)=1.$$
Hence: $1 \in B_2$ and, consequently, $2 \in B_2$. But then:
$$ y(1)=2$$
$$y(l)=l+1.$$ Therefore $(l+1) \in B_2$, which is a contradiction. Hence
$(l+2)$ does not belong to $B_2$ and the same can be said about
$(l+3),(l+4),...,l+n-1$.

On the other hand, if some $t \in \{3,...,l-1\}$ belongs to $B_2$
we have $$y(t)=t$$
$$y(l)=l+1.$$ Then, $(l+1) \in B_2$ (contradiction).

Finally, if $r \in \{1,2\}$ belongs to $B_2$:
$$y^2(r)=r$$ $$y^2(l)=l+2.$$ Hence, $(l+2) \in B_2$ (absurd, as we have seen before).

Therefore $B_2=\{l\}$. This means that all the blocks must have
only one element too because if, for instance, $a,b \in B$ and $w$
is the permutation $x^k$, for some $k \in \mathbb{N}$, that sends
$a$ to $l$, then $w(a)$ and $w(b)$ must belong to the same block.
Hence, in this case, $w(a)=l=w(b)$. That is not possible because
we are assuming that $a \neq b$.

It follows that $G=<x,y>$ is a primitive group of degree $l+n-1$.
Since $x$ is a cycle of order $l$ and $l+n-1-l+1=n$, the group $G$
is $n-transitive$ (Theorem \ref{theorem b}). For $n>5$, $G$ must
be the alternating or symmetric group (Theorem \ref{theorem d}).
If $n=4$ we have:
$$x=(1,2,...,l)$$
$$y=(1,2)(l,l+1,l+2,l+3).$$
Then, the degree of $y^2=(l,l+2)(l+1,l+3)$ is four and the group
must be the alternating or the symmetric group (Theorem
\ref{theorem a}). Since $x$ is an odd permutation, $G=S_{l+n-1}$.
Then ${\cal H}=(S_{l+n-1},x,y)$ is a hypermap with extreme duality
index because $y$ is an even permutation (see Theorem \ref{sn}).

\end{itemize}
\item[d)] If $n=l$ (and even) we can use the same generators as in
c):
$$x=(1,2,...,l)$$
$$y=(1,2)(l,l+1,l+2,...,l+n-1)=(1,2)(l,l+1,...,2l-1).$$

$$ord(x)=l$$
$$ord(y)=[2,l]=l$$
Hence $G=<x,y>=S_{l+n-1}=S_{2l-1}$. Because $x$ is an odd
permutation and $y$ is even, ${\cal H}=(S_{2l-1},x,y)$ must be a
hypermap with extreme duality index (see Theorem \ref{sn}).

\item[e)] Suppose $n=l$ (both odd).

Let the generators of the group be these two even permutations in
$S_{l+1}$:
$$x=(1,2,...,l)$$
$$y=(2,3,...,l+1).$$

The group generated by these two elements is primitive since it is
a 2-transitive group (Theorem \ref{theorem c}). Moreover,
$x^{-1}y^{-1}xy=(1,l+1)(l-1,l)$ is a permutation that just moves
four elements. It follows that $G=\langle x, y \rangle$ is a group
of minimal degree $\leq 4$. Hence, by Theorem \ref{theorem a},
$\langle x,y \rangle= A_l$ if $l+1 \geq 8$, i.e. $l \geq 7$. Cases
for low $l$ are easy to solve. $\hfill \square$

\end{itemize}
\bigskip

With this proof, we have not only shown that it is possible to
find the required hypermaps but also to describe them, since the
generators of the alternating or symmetric groups are explicitly
given for each case. We have also proved a slightly weaker result:

\begin{theorem}
For any $l$, $n \in \mathbb{N}$, one can find two permutations $x$
and $y$ (of orders $l$ and $n$, respectively), such that $x$ and
$y$ generate an alternating or a symmetric group. $\hfill \square$
\end{theorem}

\end{document}